\documentclass[11pt]{article}
\usepackage{amsmath,enumerate,fancyhdr,amssymb, amsthm, pstricks, epsf, lscape}

\newcommand{\commentout}[1]{}
\newcommand{\co}[1]{}

\co{
\setlength{\oddsidemargin}{0cm} 
\setlength{\evensidemargin}{0cm}
\setlength{\textwidth}{16.0cm} 
\setlength{\topmargin}{1cm}
\setlength{\textheight}{19cm}
}

\def\eref#1{(\ref{#1})}

\newcommand{\transpose}{\mbox{${}^{\text{T}}$}}

\newcommand{\norm}[1]{\parallel \! #1 \! \parallel}

\newcommand{\rad}[1]{\mathbb{R}^{#1}}

\newcommand{\zad}[1]{\mathbb{Z}^{#1}}

\newcommand{\latt}[1]{L(#1)}

\newcommand{\la}{\langle}
\newcommand{\ra}{\rangle}

\newcommand{\nin}{\noindent}

\newcommand{\bv}{\bar{v}}

\newcommand{\bd}{\bar{d}}
\newcommand{\bV}{\bar{V}}

\newcommand{\colmax}{\operatorname{colmax}}

\newcommand{\mydet}{\operatorname{det}}

\newcommand{\lin}{\operatorname{lin}}

\newcommand{\Proj}{\operatorname{Proj}}

\newtheorem{Definition}{Definition}

\newtheorem{Setup}{Setup}
\newtheorem{Example}{Example}
\newtheorem{Counterexample}{Counterexample}
\newtheorem{Proposition}{Proposition}
\newtheorem{Lemma}{Lemma}
\newtheorem{Theorem}{Theorem}
\newtheorem{Corollary}[Definition]{Corollary}
\newtheorem{Remark}[Definition]{Remark}
\newtheorem{Assumption}{Assumption}
\newtheorem{Recipe}{Recipe}

\newcommand{\beq}{\begin{equation}}
\newcommand{\eeq}{\end{equation}}
\newcommand{\beqa}{\begin{eqnarray}}
\newcommand{\eeqa}{\end{eqnarray}}
\newcommand{\ba}{\begin{array}}
\newcommand{\ea}{\end{array}}
\newcommand{\bac}{\begin{array}{ccccccccccc}}
\newcommand{\eac}{\end{array}}
\newcommand{\bprop}{\begin{Proposition}}
\newcommand{\eprop}{\end{Proposition}}

\newcommand{\bcex}{\begin{Counterexample}}
\newcommand{\ecex}{\end{Counterexample}}

\newcommand{\beqast}{\begin{eqnarray*}}
\newcommand{\eeqast}{\end{eqnarray*}}
\newcommand{\benum}{\begin{enumerate}}
\newcommand{\eenum}{\end{enumerate}}
\newcommand{\bit}{\begin{itemize}}
\newcommand{\eit}{\end{itemize}}
\newcommand{\bth}{\begin{Theorem}}
\newcommand{\enth}{\end{Theorem}}

\newcommand{\bdef}{\begin{Definition}}
\newcommand{\Edef}{\end{Definition}}

\newcommand{\bsetup}{\begin{Setup}}
\newcommand{\esetup}{\end{Setup}}
\newcommand{\ble}{\begin{Lemma}}
\newcommand{\ele}{\end{Lemma}}
\newcommand{\bex}{\begin{Example}}
\newcommand{\eex}{\end{Example}}
\newcommand{\bcor}{\begin{Corollary}}
\newcommand{\ecor}{\end{Corollary}}
\newcommand{\brem}{\begin{Remark}}
\newcommand{\erem}{\end{Remark}}
\newcommand{\bass}{\begin{Assumption}}
\newcommand{\eass}{\end{Assumption}}

\newcommand{\brep}{\begin{Recipe}}
\newcommand{\erep}{\end{Recipe}}

\newcommand{\pf}[1]{\vspace{.35cm} \nin {\bf Proof {#1} }}

\newcommand{\bpx}{\begin{pmatrix}}
\newcommand{\epx}{\end{pmatrix}}

\newcommand{\bbx}{\begin{bmatrix}}
\newcommand{\ebx}{\end{bmatrix}}

\alph{footnote}						

\date{\today}

\begin{document}
\pagenumbering{roman}			

\thispagestyle{empty}				

\newcommand{\Nat}{I\!\!N}
\newcommand{\RR}{I\!\!R}
\newcommand{\ZZ}{I\!\!Z}

\title{\bf On sublattice determinants in reduced bases} 
\author{G\'{a}bor Pataki and Mustafa Tural \thanks{Department of Statistics and Operations Research, UNC Chapel Hill, {\bf gabor@unc.edu, tural@email.unc.edu}} \\
Technical Report 2008-02 \\ Department of Statistics and Operations Research, UNC Chapel Hill}

\date{}

\maketitle								


\begin{abstract}
\nin Lenstra, Lenstra, and Lov\'asz in \cite{LLL82} proved several inequalities showing that the vectors in an LLL-reduced basis 
are short, and near orthogonal. 
Here we  present generalizations, from which with $k=1, \,$ and $k=n$ we can recover their inequalities:
\bth \label{sublattice-dets}
Let  $b_1,\dots,b_n \in \rad{m} $ be an LLL-reduced basis of the lattice $L, \, $ and $d_1,\dots, d_k$ arbitrary linearly independent vectors in $L$. Then 
\beqa \label{1}
\norm{b_{1}} & \leq & 2^{(n-k)/2 + (k-1)/4} (\mydet \latt{d_1, \dots, d_k})^{1/k}, \\
\label{2}
\mydet \latt{b_1, \dots, b_k} & \leq & 2^{k(n-k)/2} \mydet \latt{d_1, \dots, d_k}, \\
\label{3}
\mydet \latt{b_1, \dots, b_k}  & \leq & 2^{k(n-k)/4} (\mydet L)^{k/n}, \\
\label{4}
\norm{b_1} \cdots \norm{b_k}  & \leq  & 2^{k(n-k)/2 + k(k-1)/4} \mydet \latt{d_1, \dots, d_k}, \\
\label{5}
\norm{b_1} \cdots \norm{b_k}  & \leq  & 2^{k(n-1)/4} (\mydet L)^{k/n}.
\eeqa
\enth
\qed

\nin In the most general setting, we prove: 
\bth \label{sublattice-dets2}
Let  $b_1,\dots,b_n \in \rad{m} $ be an LLL-reduced basis of the lattice $L, \, 1 \leq k \leq j \leq n, $ and $d_1,\dots, d_j$ arbitrary linearly independent vectors in $L$. 
Then 
\beqa 
\label{6} 
\mydet \latt{b_1, \dots, b_k}  & \leq & 2^{k(n-j)/2 + k(j-k)/4} (\mydet \latt{d_1, \dots, d_j})^{k/j}, \\ 
\label{7} 
\norm{b_1} \cdots \norm{b_k}   & \leq & 2^{k(n-j)/2 + k(j-1)/4} (\mydet \latt{d_1, \dots, d_j})^{k/j}. 
\eeqa
\enth
\qed
\vspace{1cm}


\end{abstract}


\newpage								
\pagenumbering{arabic}			

\section{Lattices and Basis Reduction}
A lattice in $\rad{m}$ is a set of the form 
\beq \label{defL}
L \, = \, \latt{b_1, \dots, b_n} \, = \, \left\{ \, \sum_{i=1}^n \lambda_i b_i \, | \, \lambda_i \in \zad{}, \, (i=1, \dots, m) \, \right\},
\eeq
where $b_1, \dots, b_n$ are linearly independent vectors in $\rad{m}, $ and are called a {\em basis} of $L$. 
If $B = [b_1, \dots, b_n], \,$ then we also call $B$ a basis of $L$, and write $L = \latt{B}.$ The determinant of $L$ is 
\beq \label{det}	
\mydet L = \sqrt{\mydet B \transpose B},
\eeq
where $B$ is a basis of $L, \,$ with $\mydet L \,$ actually independent of the choice of $B$. 

Finding a short, nonzero vector in a lattice is a fundamental algorithmic problem with many uses in cryptography, optimization, and number theory. 
For surveys we refer to \cite{GLS93}, \cite{K87agn}, \cite{S86}, and \cite{MI02}. 
More generally, one may want to find a reduced basis consisting of short, and nearly orthogonal vectors.

A basis $b_1, \dots, b_n$ that is reduced according to the definition of 
Lenstra, Lenstra, and Lov\'asz \cite{LLL82} is computable in polynomial time in the case of rational lattices, and the $b_i$ are reasonably short, and near orthogonal, namely
 \beqa  \label{lll1} 
\norm{b_1} & \leq & 2^{(n-1)/4} (\mydet L)^{1/n},  \,  \\ 
\label{lll2} 
\norm{b_1} & \leq &  2^{(n-1)/2} \norm{d} \, \text{for \, any \,} d \in L \setminus \{ 0 \}, \\ 
\label{lll3} 
\norm{b_1} \cdots \norm{b_n}   & \leq & 2^{n(n-1)/4} \mydet L.
\eeqa
hold. Korkhine-Zolotarev (KZ) bases, which were described in \cite{KZ1873} by Korkhine, and Zolotarev, and by Kannan in \cite{K87}
have stronger reducedness properties,  for instance, the first vector in a KZ basis is the shortest vector of the lattice.
However, KZ bases are computable in polynomial time only when $n$ is fixed. 
Block KZ bases proposed by Schnorr in \cite{Sh87} form a hierarchy in between: one can trade on the quality of the basis to 
gain faster computing times. 

Our Theorem 1 generalizes inequalities \eref{lll1} through \eref{lll3}. For instance, \eref{1} with $k=n \,$ yields \eref{lll1}, and with 
$k=1 \,$ yields \eref{lll2}. In turn, from \eref{6} in Theorem 2 with $j=k, \,$ and from \eref{7} with $j=n \,$ we recover the inequalities of Theorem 1.

It would be interesting to see whether stronger versions of our results can be stated for KZ, or block KZ bases.

As a tool we use Lemma 1 below,  which may be of independent interest. For $k=1$ we can recover  from it 
Lemma (5.3.11) in \cite{GLS93} (proven as part of Proposition (1.11) in \cite{LLL82}).
To state it, we will recall the notion of Gram-Schmidt orthogonalization.
If $b_1, \dots, b_n \in \rad{m}$ is a basis of $L$, then the corresponding 
Gram-Schmidt vectors $b_1^*, \dots, b_n^*, $ are defined as 
\beq \label{gso} 
b_1^*=b_1 \; \text{and} \;
b_i^*=b_i-\sum_{j=1}^{i-1} \mu_{ij} b_j \; \text{for} \; i=1, \dots, n-1, 
\eeq
with $\mu_{ij}=  \la b_i, b_j^* \ra/ \la b_j^*, b_j^* \ra, $ where $\la. , . \ra $ is the usual inner product on $\rad{m}$ .

\ble \label{index-lemma=det}
Let $d_1, \dots, d_k$ be linearly independent vectors from the lattice   $L,$
and $b_1^*, \dots, b_n^*$ the Gram Schmidt orthogonalization of an arbitary basis. 
Then
\beqa
\mydet \latt{d_1, \dots, d_k} \geq \min_{1 \leq i_1 <  \dots < i_k \leq n} \left\lbrace  \norm{b_{i_1}^*} \dots \norm{b_{i_k}^*} \right\rbrace.
\eeqa
\ele
\qed

In the rest of this section we collect necessary definitions, and results. In Section \ref{sect-inter} we prove Lemma 1, and in Section \ref{sect-thm1} we prove 
Theorem \ref{sublattice-dets2}. 

We call $b_1, \dots, b_n$ an {\em LLL-reduced basis of } $L,$ if 
 \beqa \label{mucond}
| \mu_{ji} | & \leq & 1/2 \; \,\,\, (j=2, \dots, n; \, i = 1, \dots, j-1), \, \text{and} \\ \label{exch-cond}
\norm{b_j^*+ \mu_{j,j-1}b_{j-1}^*}^2 & \geq &  3/4 \norm{b_{j-1}^*}^2 \,\, (1<j \leq n).
\eeqa
From (\ref{mucond}) and (\ref{exch-cond}) it follows that 
\beq \label{new-cond2}
\norm{b_{i}^*}^2  \; \leq \;  2^{j-i} \norm{b_{j}^*}^2 \; (1 \leq i \leq j \leq n).
\eeq
If $b_1, \dots, b_n$ are linearly independent vectors, then 
\beq \label{bndet}
\mydet L(b_1, \dots, b_n) \, =  \, \mydet L(b_1, \dots, b_{n-1}) \norm{b'}, 
\eeq
where $b'$ is the projection of $b_n$ on the orthogonal complement of the linear span of $b_1, \dots, b_{n-1}$. 

An integral square matrix $U$ with $\pm 1$ determinant is called unimodular. 
An elementary column operation performed on a matrix $A$ is either 1) exchanging two columns, 2) multiplying a column by $-1$, or 3) adding an integral multiple 
of a column to another column. Multiplying a matrix $A$ from the right by a unimodular $U$ is equivalent to performing a sequence of elementary column operations on $A$.

\section{Proof of Lemma 1} 
\label{sect-inter}

We need the following 

\nin{\bf Claim}  
There are elementary column operations performed on $d_1, \dots, d_k \,$ that yield $\bd_1, \dots, \bd_k$ with 
\beqa \label{bardi} 
\bd_i=\sum_{j=1}^{t_i} \lambda_{ij} b_j \; \text{for} \; i=1,\dots,k,
\eeqa
where $\lambda_{ij} \in \zad{}, \, \lambda_{i, t_i} \neq 0, $ and 
\beq \label{tk} 
t_k > t_{k-1} > \dots > t_1. 
\eeq

\pf{of Claim} Let us write 
\beqa
B V & = & [ d_1, \dots, d_k],
\eeqa
with $V$  an integral matrix. 
Analogously to how the Hermite Normal Form of an integral matrix is computed, we can 
do elementary column operations on $V$ to obtain $\bV$ with 
\beq
t_k := \max \, \{ \, i \, | \, \bv_{ik} \neq 0 \, \} \, > \, t_{k-1} := \max \, \{ \, i \, | \, \bv_{i,k-1} \neq 0 \, \} \, > \, \dots \, > \, t_1 := \max \, \{ \, i \, | \, \bv_{i1} \neq 0 \, \}.
\eeq
Performing the same elementary column operations on $d_1, \dots, d_k$ yield $\bd_1, \dots, \bd_k \,$  which satisfy 
\beqa
B \bV & = & [ \bd_1, \dots, \bd_k],
\eeqa
so they satisfy \eref{bardi}.

\nin{\bf End of proof of Claim} 

\nin Obviously 
\beq \label{deteq}
\mydet \, \, \latt{ \bar{d}_1, \dots, \bar{d}_k } \, =  \, \mydet \, \, \latt{ d_1, \dots, d_k }.
\eeq
Substituting from \eref{gso} for $b_i$  we can rewrite \eref{bardi} as 
\beqa \label{bardi2} 
\bar{d}_i=\sum_{j=1}^{t_i} \lambda_{ij}^* b_j^* \; \text{for} \; i=1,\dots,k,
\eeqa
where the $\lambda_{ij}^*$ are now reals, but $\lambda_{i, t_i}^* = \lambda_{i, t_i}$ nonzero integers. 

\nin For all $i$ we have 
\beq
\lin \, \{ \, \bd_1, \dots, \bd_{i-1} \, \} \, \subseteq \, \lin \{ \, b_1^*, \dots, b_{t_{i-1}}^* \, \}.
\eeq
Therefore
\beq \label{bdi} 
\norm{\Proj{ \, \{ \, \bd_i \, | \, \{ \, \bd_1, \dots, \bd_{i-1} \, \}^\perp \, \}}} \, \geq \, \norm{\Proj{ \, \{ \, \bd_i \, | \, \{ \, b_1^*, \dots, b_{t_{i-1}}^* \, \}^\perp \, \}}} \, \geq \, \norm{\lambda_{i, t_i} b_{t_i}^*} \, \geq \, \norm{b_{t_i}^*}
\eeq
holds, with the second inequality coming from  \eref{tk}. So applying \eref{bndet} repeatedly we get
\beq 
\ba{rcl}
\mydet \, \, \latt{ \bar{d}_1, \dots, \bar{d}_k } & \geq & \mydet \latt{\bar{d}_1, \dots, \bar{d}_{k-1}} \norm{b_{t_k}^*} \\
                                          & \dots &  \\
                                          & \geq &   \norm{b_{t_1}^*} \norm{b_{t_2}^*} \dots \norm{b_{t_k}^*},
\ea
\eeq
which together with \eref{deteq} completes the proof.
\qed

\section{Proof of Theorem 1 and Theorem 2} 
\label{sect-thm1} 

The plan of the proof is as follows: we first prove \eref{1} through \eref{3} in Theorem 1. Then we prove 
Theorem 2. Finally, \eref{4} follows as a special case of \eref{7} with $j=k$; and \eref{5} as a special case of \eref{7} with $j=n$. 

\pf{of \eref{1} and \eref{2}} 
Lemma  \ref{index-lemma=det} implies 
\beqa
\mydet \, \, \latt{ d_1, \dots, d_k } & \geq & \norm{b_{t_1}^*} \norm{b_{t_2}^*} \dots \norm{b_{t_k}^*}
\eeqa
for some $t_1, \dots, t_k \in \{1, \dots, n \}$ distinct indices. 
Clearly 
\beq \label{t1tk}
t_1 + \dots + t_k \leq kn - k(k-1)/2
\eeq
holds.
Applying first \eref{new-cond2}, then \eref{t1tk} yields 
\beq
\ba{rcl}
(\mydet \, \, \latt{ d_1, \dots, d_k })^2 & \geq & \norm{b_{1}^*}^2  2^{(1 - t_1)} \dots \norm{b_{1}^*}^2 2^{(1 - t_k)} \\
                                      & = &  \norm{b_{1}^*}^{2k} 2^{k - (t_1 + \dots + t_k)} \\ 
       & \geq & \norm{b_{1}}^{2k}  2^{k(k+1)/2 - kn},
\ea
\eeq
which is equivalent to \eref{1}. 
Similarly, 
\beq
\ba{rcl}
(\mydet \, \, \latt{ d_1, \dots, d_k })^2 & \geq & \norm{b_{1}^*}^2  2^{(1 - t_1)} \norm{b_{2}^*}^2  2^{(2 - t_2)}  \dots \norm{b_{k}^*}^2 2^{(k - t_k)} \\
                                      & = &  \norm{b_{1}^*}^2 \dots \norm{b_{k}^*}^2 2^{(1+ \dots + k) - (t_1 + \dots + t_k)} \\
       & \geq & \norm{b_{1}^*}^2 \dots \norm{b_{k}^*}^2   2^{k(k-n)},
\ea
\eeq
which is equivalent to \eref{2}. 

\qed

\pf{of \eref{3}} The proof is  by induction. Let us write $D_k = (\mydet \latt{b_1, \dots, b_k})^2$.  For $k=n-1, \,$ multiplying the inequalities
\beq
\norm{b_{i}^*}^2 \, \leq \, 2^{n-i} \norm{b_{n}^*}^2 \, (\; i=1,\dots,n-1)
\eeq
gives 
\beqa
D_{n-1} & \leq & 2^{n(n-1)/2} (\norm{b_{n}^*}^{2})^{n-1} \\
      &  =   & 2^{n(n-1)/2} \left( \dfrac{D_n}{D_{n-1}} \right)^{n-1},
\eeqa
and after simplifying, we get 
\beqa \label{dn-1}
D_{n-1} & \leq & 2^{(n-1)/2} (D_n)^{1 - 1/n}.
\eeqa

Suppose that  \eref{3} is true for $k \leq n-1$; we will prove it for $k -1$. 
Since $b_1, \dots, b_{k}$ forms an LLL-reduced basis of $\latt{b_1,\dots, b_k}$ we can replace $n$ by $k$ in (\ref{dn-1}) to get 
\beqa \label{dl-1}
D_{k-1} & \leq & 2^{(k-1)/2} (D_k)^{(k-1)/k}.
\eeqa
By the induction hypothesis, 
\beqa \label{dldn-2}
D_k & \leq & 2^{k(n-k)/2} (D_n)^{k/n},
\eeqa
from which we obtain
\beqa \label{dldn-3}
(D_k)^{(k-1)/k} & \leq & 2^{(k-1)(n-k)/2} (D_n)^{(k-1)/n}.
\eeqa
Using the upper bound on $(D_k)^{(k-1)/k}$ from (\ref{dldn-3}) in (\ref{dl-1}) yields 
\beqa \label{dl-1-2}
D_{k-1} & \leq & 2^{(k-1)/2}  2^{(k-1)(n-k)/2} (D_n)^{(k-1)/k} \\
             & =    & 2^{(k-1)(n-(k-1))/2} (D_n)^{(k-1)/n},
\eeqa
as required. 

\qed

\pf{of Theorem 2}  From \eref{3} and \eref{2} we have 
\beqa
\label{one} \mydet \latt{b_1, \dots, b_k} & \leq & 2^{k(j-k)/4} (\mydet \latt{b_1, \dots, b_j})^{k/j}, \\
\label{two} \mydet \latt{b_1, \dots, b_j} & \leq & 2^{j(n-j)/2} \mydet \latt{d_1, \dots, d_j}.
\eeqa
Raising \eref{two} to the power of $k/j \,$ gives
\beqa
\label{three} (\mydet \latt{b_1, \dots, b_j})^{k/j} & \leq & 2^{k(n-j)/2} \mydet (\latt{d_1, \dots, d_j})^{k/j},  
\eeqa
and plugging \eref{three} into \eref{one} proves \eref{6}. 

It is shown in \cite{LLL82} that 
\beqa \label{bibd} 
\norm{b_i}^2 & \leq & 2^{i-1} \norm{b_i^*}^2 \; \text{for \;} i=1, \dots, n.
\eeqa
Multiplying these inequalities for $i=1, \dots, k \,$ yields
\beqa \label{bibd2} 
\norm{b_1} \cdots \norm{b_k}  & \leq  & 2^{k(n-1)/4} \mydet \latt{b_1, \dots, b_k},
\eeqa
and using \eref{bibd2} with \eref{6} yields \eref{7}. 

\qed

\brem {\rm
The $k$th successive minimum of $L$ is defined as the smallest real number $t$, such that 
there are $k$ linearly independent vectors in $L$ with length bounded by $t$. It is denoted by $\lambda_k(L)$.
With the same setup as for \eref{lll1}-\eref{lll3} it is shown in \cite{LLL82} that 
\beqa \label{succ} 
\norm{b_i} & \leq & 2^{n-1} \lambda_i(L) \; \text{for} \, i=1, \dots, n.
\eeqa
For KZ, and block KZ bases similar results were shown in \cite{LLS90},  and \cite{Sh94}, resp.

The successive minimum results \eref{succ}  give a more global, and refined  view of the lattice, and the reduced basis, 
than \eref{lll1} through \eref{lll3}. Our Theorems 1 and 2 are similar in this respect, but they seem to be independent of \eref{succ}. Of course, multiplying 
the latter for $i=1, \dots, k \,$ gives an upper bound on $\norm{b_1} \cdots \norm{b_k}, \,$ but in different terms. 

The quantites $\mydet \latt{b_1, \dots, b_k}$ and $\norm{b_1} \dots \norm{b_k}$ are also connected by 
\beqa \label{normp}
\mydet \latt{b_1, \dots, b_k} & = & \norm{b_1} \dots \norm{b_k} \sin \theta_2 \dots \sin \theta_k,
\eeqa
where $\theta_i$ is the angle of $b_i$ with the subspace spanned by $b_1, \dots, b_{i-1}$. 
In \cite{Ba86} Babai showed that the sine of the angle of {\em any} basis vector with the subspace spanned by the other basis vectors in a $d$-dimensional lattice 
is at least $(\sqrt{2}/3)^d$. One could combine the lower bounds on $\sin \theta_i$ with the upper bounds on 
$\mydet \latt{b_1, \dots, b_k} \,$ to find an upper bound on $\norm{b_1} \dots \norm{b_k}. \,$ However, the result would be weaker than 
\eref{4} and \eref{5}. 
}
\erem

\nin{\bf Acknowledgement} The first  author thanks Ravi Kannan for helpful discussions.

\bibliography{IP_Refs}

\end{document}